\documentclass[10pt,draftcls,onecolumn,journal,transmag]{IEEEtran}

\usepackage{lipsum}
\usepackage{graphicx}
\ifCLASSOPTIONcompsoc
\usepackage[caption=false, font=normalsize, labelfont=sf, textfont=sf]{subfig}
\else
\usepackage[caption=false, font=footnotesize]{subfig}
\fi

\ifCLASSINFOpdf
\else
\fi
\usepackage{amsmath}

\usepackage{algorithmic}

\usepackage{amssymb}
\usepackage{graphicx}
\usepackage{graphicx}
\usepackage{graphicx}
\usepackage{amsmath}
\usepackage{epstopdf} 

\usepackage{multirow}
\usepackage{mathptmx}      
\usepackage{amsthm}
\usepackage{geometry}
\usepackage{lineno}
\usepackage{color}
\usepackage{float}

\newtheorem{definition}{Definition}
\newtheorem{theorem}{Theorem}
\newtheorem{example}{Example}

\newtheorem{proposition}{Proposition}
\newtheorem{remark}{Remark}

\newtheorem{corollary}{corollary}

\begin{document}
%
\title{\huge{Piecewise Sparse Recovery in Unions of Bases}}


\author{\IEEEauthorblockN{Chong-Jun Li\IEEEauthorrefmark{1}and
Yi-Jun Zhong}\\
{\IEEEauthorrefmark{1}School of Mathematical Sciences,
Dalian University of Technology, Dalian 116024, China}
\thanks{This work was supported by the National Natural Science Foundation of China (Grant Nos. 11871137, 11572081), and the Fundamental Research Funds for the Central Universities of China (Grant No. QYWKC2018007).\quad
Corresponding author: Chong-Jun Li (Email:  chongjun@dlut.edu.cn).}}

\date{ }

%



\maketitle

\begin{abstract}
Sparse recovery is widely applied in many fields, since many signals or vectors can be sparsely represented under some frames or dictionaries. Most of fast algorithms at present are based on
solving $l^0$ or $l^1$ minimization problems and they are efficient in sparse recovery. However, compared with the practical results, the theoretical sufficient conditions on the sparsity of the signal for
$l^0$ or $l^1$ minimization problems and algorithms are too strict. \par

In many applications, there are signals with certain structures as piecewise sparsity.
Piecewise sparsity means that the sparse signal $\mathbf{x}$ is a union of several sparse sub-signals, i.e., $\mathbf{x}=(\mathbf{x}_1^T,\ldots,\mathbf{x}_N^T)^T$, corresponding to the matrix $A$ which is composed of union of bases $A=[A_1,\ldots,A_N]$. In this paper, we consider the uniqueness and feasible conditions for piecewise sparse recovery. We introduce the mutual coherence for the sub-matrices $A_i\ (i=1,\ldots,N)$ to study the new upper bounds of $\|\mathbf{x}\|_0$ (number of nonzero entries of signal) recovered by $l^0$ or $l^1$ optimizations.
The structured information of measurement matrix $A$ is used to improve the sufficient conditions for successful piecewise sparse recovery and also improve the reliability of $l_0$ and $l_1$ optimization models on recovering global sparse vectors.
\end{abstract}

\begin{center}
\end{center}

\begin{center}
{\bf MR(2010) Subject Classification}\\ 90C25; 94A12
\end{center}

%
\IEEEpeerreviewmaketitle

\section{Introduction}
%
%
%
%
\IEEEPARstart{I}{n} this paper,  we consider recovering a sparse signal $\mathbf{x}^*\in \mathbf{R}^n$ from an underdetermined system of linear equation
\begin{equation}\label{originalmodel}
  A\mathbf{x}^*=\mathbf{b},
\end{equation}
where $\mathbf{b}\in\mathbf{R}^m$ is a measurement vector, $A\in \mathbf{R}^{m\times n}$ is a measurement matrix. If the vector $\mathbf{x}^*$ has at most $s\leq m < n$ nonzero entries, then it is named as $s$-sparse vector, the corresponding index set of nonzero entries is called support $\mathbf{S}={\rm supp}(\mathbf{x}^*)$.
There are many theories, algorithms and applications on this problem of sparse recovery \cite{Elad2010Book}.

One approach to find the sparsest solution of Eq.~(\ref{originalmodel}) is greedy algorithm (GA), which solves the following $l^0$ minimizing
solution, named as P$_0$ problem:
\begin{equation}\label{l0model}
  \min_{\mathbf{x}}\|\mathbf{x}\|_0,~~~~~s.t.~~ A\mathbf{x}=\mathbf{b}.
\end{equation}
One of the most popular greedy methods is the orthogonal matching
pursuit (OMP) as proposed in \cite{mallat1,pati,tropp2}. It iteratively adds components to the support
of the approximation $\mathbf{x}^k$ whose correlation to the current residual is maximal. There are many other greedy methods for sparse
recovery, for example, iterative hard thresholding (IHT) \cite{blumensath},
stagewise OMP (StOMP) \cite{Donoho2012}, regularized OMP (ROMP) \cite{Needell2008-1,Needell2010},
compressive sampling matching pursuit (CoSaMP) \cite{Needell2008}, subspace pursuit (SP) \cite{dai},
iterative thresholding with inversion (ITI) \cite{maleki}, hard thresholding pursuit (HTP)
\cite{foucart} etc.

Another approach is convex relaxation which
solves a convex program whose minimizer is obtained to approximate the target signal. The basis pursuit gains lots of attention which determines the sparsest representation of $\mathbf{x}^*$ by solving the following $l^1$ minimization problem, named as P$_1$ problem or BP problem (method):
\begin{equation}\label{l1model}
  \min_{\mathbf{x}}\|\mathbf{x}\|_1,~~~~~s.t.~~ A\mathbf{x}=\mathbf{b}.
\end{equation}
 Many algorithms have been proposed to complete the optimization, including interior-point methods \cite{Kim2007}, projected gradient methods \cite{Figueiredo2007}, and iterative thresholding \cite{Daubechies2004} etc. \par

There are three fundamental problems concerned in this paper:
\begin{enumerate}
\item Uniqueness of solution of the P$_0$ problem.
\item Feasibility of GA (or OMP) for solving the P$_0$ problem.
\item Equivalence between the P$_1$ problem and the P$_0$ problem, or
   feasibility of BP method to obtain the sparsest solution.
\end{enumerate}

There are several tools proposed to formalized the notion of the suitability for the above three problems, such as the mutual coherence \cite{Donoho2001}, the spark \cite{Donoho2003}, the cumulative coherence \cite{Tropp2004}, the exact recovery condition (ERC) \cite{Tropp2004}, and the restricted isometry constants (RICs) \cite{Candes2005,Candes2006-2,Candes2008-1}. It is well-known that the sufficient and necessary condition for the uniqueness of the solution of the P$_0$ problem (\ref{l0model}) is (\cite{Donoho2003})
\begin{equation}\label{sparkcondition}
\|\mathbf{x}\|_0<{\rm spark}(A)/2,
\end{equation}
or the RIC of the matrix $A$ satisfies
$\delta_{2s}<1$ and $\|\mathbf{x}\|_0\leq s$ (\cite{Candes2005}). The equivalence between the P$_1$ model and the P$_0$ model is guaranteed by
$\delta_{2s}<\sqrt{2}-1$ (\cite{Candes2006-2,Candes2008-1}). For a given matrix or dictionary $A$, however, it is difficult to compute the spark or verify the RIP conditions. By contrast, we can easily compute the mutual coherence of matrix.
The general case discussed in \cite{Donoho2001} showed
that one sufficient condition for the uniqueness of the solution of the P$_0$ problem (\ref{l0model}) is
\begin{equation}\label{strictbound}
\|\mathbf{x}\|_0<\frac{1}{2}(1+\frac{1}{\mu}),
\end{equation}
where $\mu$ is the mutual coherence of measurement matrix $A$. Furthermore, the condition Eq. (\ref{strictbound}) is also a sufficient condition which ensures the OMP (greedy method) and BP methods for recovering the optimal $s$-sparse solution \cite{Tropp2004}. However, in applications, OMP or BP method can work well even when the condition (\ref{strictbound}) is not satisfied, i.e., when
\begin{equation*}
 \frac{1}{2}(1+1/\mu(A))\leq \|\mathbf{x}\|_0<{\rm spark}(A)/2,
\end{equation*}
which means the sufficient condition (\ref{strictbound}) is strict for sparse recovery to some extent, or the "gap" between the optimal upper bound $\frac{1}{2}{\rm spark}(A)$ and the practical upper bound $\frac{1}{2}(1+1/\mu(A))$ is big.\par
To our surprise, the result in \cite{elad2002} shows that uniqueness of the $l^0$ minimization P$_0$ problem solution
can be achieved for improved condition
$\|\mathbf{x}\|_0<\frac{1}{\mu}$.
The authors also showed that the solutions of the P$_0$ and P$_1$ problems
coincide for $\|\mathbf{x}\|_0<\frac{\sqrt{2}-0.5}{\mu}$. These two improved conditions were obtained in the special case where $A$ is in pairs of orthogonal
bases.
It was also shown in \cite{gribonval2003,Tropp2004} that if the matrix $A$ is a union of $N(\geq 2)$ orthogonal bases, improved conditions are possible. The sufficient condition for OMP to solving the P$_0$ model for $N$ orthogonal bases was improved with
\begin{equation}\label{orthouni}
   \|\mathbf{x}\|_0<(\frac{1}{2}+\frac{1}{2(N-1)})\frac{1}{\mu},
\end{equation}
and for BP to solving P$_1$ problem was improved with
\begin{equation}\label{Gribonvalresult}
  \|\mathbf{x}\|_0<(\sqrt{2}-1+\frac{1}{2(N-1)})\frac{1}{\mu}.
\end{equation}
Especially, when the vector $\mathbf{x}$ is $(s_1,\ldots,s_N)$-piecewise sparse (see Definition \ref{Def1}, $\|\mathbf{x}\|_0=s_1+\cdots+s_N$) or vector $\mathbf{b}$ is a superposition of $s_i$ atoms from the $i$-th basis,
  and $s_1\leq \cdots \leq s_N$ with $A$ in union of $N$ orthogonal bases, the exact recovery condition (ERC)
is guaranteed by
  \begin{equation}\label{Troppresult}
    \sum \limits_{j=2}^N \frac{\mu s_j}{1+\mu s_j}<\frac{1}{2(1+\mu s_1)}.
  \end{equation}
Thus, both the OMP and BP can find the sparest solution under condition (\ref{Troppresult}) \cite{Tropp2004}. When $n=2$, Eq.~(\ref{Troppresult}) is $2\mu s_1s_2+\mu s_2<1~(s_1\leq s2)$.

\begin{table}[htbp]
\begin{center}
\caption{A list of theoretical upper bounds for sparse recovery. }\label{Tab1}
\small
\begin{tabular}{ccc}
  \hline
condition index &structure of $A$ &the upper bound of $\|\mathbf{x}\|_0(=s_1+s_2)$\\
\hline
 condition 1& general case & $\|\mathbf{x}\|_0<\frac{1}{2}(1+\frac{1}{\mu(A)})$ \\
  condition 2&pair of orthogonal bases~(uniqueness) & $\|\mathbf{x}\|_0<\frac{1}{\mu(A)}$ \\
   condition 3&pair of orthogonal bases~(equivalence)&  $\|\mathbf{x}\|_0<\frac{\sqrt{2}-0.5}{\mu(A)}$\\
  condition 4&pair of orthogonal bases~(ERC)  & $2\mu^2 s_1s_2+\mu s_2<1$\\
  \hline
\end{tabular}
\end{center}
\label{conditionlist}
\end{table}
\begin{figure}
  \centering
  \includegraphics[width=7.5cm]{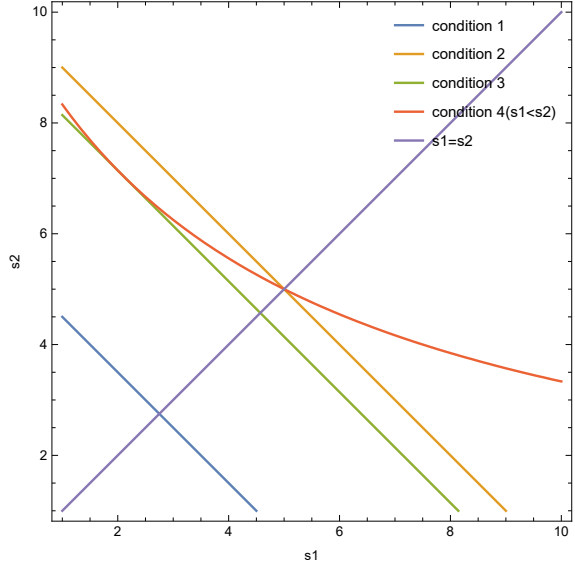}\\
  \caption{A plot of the upper bounds for the sparse recovery in Table \ref{conditionlist} with $\mu(A)=0.05$.}\label{conditionplot}
\end{figure}
The four conditions are concluded in Table \ref{Tab1} for the case of $A$ in pairs of orthogonal bases compared with the general case of $A$.
It is observed from the Fig.~\ref{conditionplot} (presented in \cite{elad2002}) that the sufficient conditions based on the mutual coherence can be improved by considering the structure of matrix $A$.\par

Note that, many natural and useful redundant dictionaries (measurement matrix $A$) cannot be written as a union of orthogonal bases.
Thus, it is necessary to study the sufficient conditions for successful recovery when
the measurement matrix $A$ in general setting, i.e, $A$ is a union of non-orthogonal bases. In this case, the corresponding vector $\mathbf{x}$ is partitioned into several parts according to the structure of matrix $A$. By considering the different sparsity of $\mathbf{x}$ according to the structure of $A$ makes it possible to deeply study the relaxed sufficient conditions for successful sparse recovery.\par

For a given vector $\mathbf{x}=(\underbrace{x_1,\ldots,x_{d_1}}_{\mathbf{x}_1^T}, \underbrace{x_{d_1+1},\ldots,x_{d_1+d_2}}_{\mathbf{x}_2^T},\ldots,\underbrace{x_{n-d_N+1},\ldots,x_{n}}_{\mathbf{x}_N^T})^T$, where
$n=\sum\limits_{i=1}^N d_i$ and $s_i=\|\mathbf{x}_i\|_0$,~$i=1,\ldots,N$.
There are three types of sparsity of vector $\mathbf{x}$:
\begin{enumerate}
\item global sparsity:~$\mathbf{x}$ is assumed to have $s=\|\mathbf{x}\|_0=\sum\limits_{i=1}^N\|\mathbf{x}_i\|_0$ nonzero entries.\par
\item block sparsity (\cite{Peotta2007,Eldar2009,Eldar2010,Elhamifar2011}):~A block
$s$-sparse vector $\mathbf{x}$ is assumed to have at
most $s$ blocks with nonzero entries, i.e., the block $l^0$ or $l^1$ norm $\|\mathbf{x}\|_{2,0}=\sum\limits_{i=1}^NI(\|\mathbf{x}_i\|_2)$ or $\|\mathbf{x}\|_{2,1}=\sum\limits_{i=1}^N\|\mathbf{x}_i\|_2$ is minimized to recovery.\par
\item \textbf{piecewise sparsity}:~as the following definition.
\end{enumerate}
\begin{definition}\label{Def1}
A vector $\mathbf{x}=(\mathbf{x}_1^T,\ldots,\mathbf{x}_N^T)^T\in \mathbf{R}^{n}$ is
partitioned into $N$ components and it is assumed that every
$\mathbf{x}_i^T\in \mathbf{R}^{n_i}$ containing nonzero entries is
sparse, $s_i=\|\mathbf{x}_i\|_0$,~$i=1,\ldots,N$. We call the vector $\mathbf{x}$ is $(s_1,\ldots,s_N)$-piecewise sparse.
\end{definition}
\par
Piecewise sparse vector is a type of vector which each part of the vector is sparse. The piecewise sparsity is different from the block sparsity. Piecewise sparse recovery are common in applications, such as
the problem of the decomposition of texture part and cartoon part of image in \cite{Starck2005}, i.e., $\mathbf{b}=A_n\mathbf{x}_n+A_t\mathbf{x}_t$ where $n$ and $t$ represent the cartoon and texture. It is assumed that both parts can be represented in some given dictionaries, thus $\mathbf{x}_n$ and $\mathbf{x}_t$ are two sparse vectors. The coefficient vector $\mathbf{x}=(\mathbf{x}_n^T,\mathbf{x}_t^T)^T$ is ``piecewise" sparse vector. Another example is
 the problem of reconstructing a surface from scattered data in
approximation space $H=\bigcup_{i=1}^N H_j$,~ where $H_j\subseteq H_{j+1}$ are
principal shift invariant (PSI)
spaces generated by a single compactly supported function \cite{hlw}, the fitting surface is $g=\sum_{i=1}^N g_i,~g_i\in H_i$ with $g_i=\sum_{j=1}^{n_i}c_j^i\phi_j^i$. The coefficients $\mathbf{c}=(\mathbf{c}^1,\mathbf{c}^2,\ldots,\mathbf{c}^N)^T$ (by $N$ pieces $\mathbf{c}^i=(c^i_1,\ldots,c^i_{n_i})^T$)
 is the vector to be determined. Due to the property of PSI spaces, the coefficients to be determined by $l^1$ minimization in \cite{hlw} are ``piecewise" sparse structured, i.e. each $\mathbf{c}^i\in \mathbf{R}^{n_i}$ is a sparse vector in $H_i$. In \cite{Zhong2018}, we firstly try to recover the piecewise sparse vector by the piecewise inverse scale space algorithm with deletion rule.  \par

It is obvious that piecewise sparsity is more general in applications, since the nonzero entries can appear in scattered
way. The corresponding matrix can be structured in a union of some
bases (orthogonal bases is a special case) $A=[A_1,\ldots,A_N]$. In this paper,
We use the mutual coherence and cumulative mutual coherence to give the conditions of piecewise sparse recovery, which
can be efficiently calculated for an arbitrary given matrix $A$ and $A_i$ ($i=1,\ldots,N$).
Inspired by the works in \cite{gribonval2003,Tropp2004}, which provide improved sufficient conditions for having unique sparse representation of signals in unions of orthogonal bases, we study the generalization of the sufficient conditions for having unique sparse representation of signals in unions of general bases corresponding to piecewise sparsity. \par


\section{Preliminaries}
We use $\mathbf{x}$ to represent a vector and $x$ represents a scalar. Define $\langle \mathbf{x},\mathbf{y}\rangle=\mathbf{x}^T\mathbf{y}$.
Let $\mathbf{x}=(\mathbf{x}_1^T,\ldots,\mathbf{x}_N^T)^T$ correspond to $A=[A_1,\ldots,A_N]$.
For convenience, let $\mathbf{S}={\rm supp}(\mathbf{x})$ and $T$ be its complement, i.e. $T = \{i:x_i=0\}$. $A_{\mathbf{S}}$ denotes the
submatrix of $A$ formed by the columns of $A$ in $\mathbf{S}$. Similarly
define $A_T$ so that $[A_{\mathbf{S}}~~ A_T]=A$.
Denote $s=|\mathbf{S}|=|supp(\mathbf{x})|$, $s_i=|S_i|=|supp(\mathbf{x}_i)|$.

\subsection{Tools used in sparse recovery}
In this part, we introduce the widely used tool for sparse
signal recovery: the \textbf{mutual coherence} of a dictionary $A\in \mathbb{R}^{m \times n}$.
Denote~$\mathbf{a}_k^i$ by the $k$-th column in the submatrix $A_i\in \mathbb{R}^{m \times n_i}$, the matrix $A$ is assumed to have unit $l^2$ norm for each column, i.e., $\|\mathbf{a}_k^i\|_2=1,\ i=1,\ldots,N,\ k=1,\ldots,n_i$.
\begin{definition}\cite{Donoho2001}
\textbf{mutual coherence of the dictionary}:
\begin{equation*}
\mu:=\max_{i,j}|\langle \mathbf{a}_{i}, \mathbf{a}_{j}\rangle|.
\end{equation*}
\end{definition}
Roughly speaking, the coherence measures how much two vectors in the dictionary can look alike. It is obvious that every orthogonal basis has coherence $\mu=0$.  A union of two orthogonal bases has coherence $\mu\geq m^{-1/2}$ \cite{Donoho2001}.
\begin{definition}\cite{Tropp2004}
the \textbf{cumulative mutual coherence} (Babel function) is defined by
\begin{equation*}
 \mu_1(s)=\max\limits_{|S|=s}\max\limits_{\mathbf{a}_i\in S^c}\sum\limits_{\mathbf{a}_j\in S}|\langle \mathbf{a}_i,\mathbf{a}_j\rangle|.
\end{equation*}
A close examination of the formula shows that $\mu_1(1)=\mu$ and that $\mu_1$ is a non--decreasing function of $s$.
\end{definition}
\begin{proposition}\cite{Tropp2004}
If a dictionary $A$ has coherence $\mu$, then $\mu_1(m)\leq m\mu$.
\end{proposition}

\begin{theorem}\cite{Tropp2004}
\textbf{Exact Recovery Condition (ERC)}. A sufficient condition for both Orthogonal Matching Pursuit and Basis Pursuit to recover the $s$-sparse $\mathbf{x}$ with $\mathbf{S}=supp(\mathbf{x})$
successfully is that
\begin{eqnarray*}
\max\limits_{j\in \mathbf{S}^C}\|(A_\mathbf{S}^TA_\mathbf{S})^{-1}A_\mathbf{S}^TA_j\|_1<1.
\end{eqnarray*}
\end{theorem}

\begin{definition}\cite{Donoho2003}
The \textbf{spark} of a matrix (dictionary) $A$ counts the least number of columns which form a linearly dependent set.
\begin{eqnarray*}
\mbox{spark}(A)=\min\limits_{\mathbf{x}\in \mbox{Ker}(A),~\mathbf{x}\neq 0}\|\mathbf{x}\|_0,
\end{eqnarray*}
where $\mbox{Ker}(A)$ is the kernel of the dictionary defined as $\mbox{Ker}(A)=\{\mathbf{x}: A\mathbf{x}=0\}$.
\end{definition}

\subsection{Tools used in piecewise sparse recovery}
Assume $A=[A_1,\ldots,A_N]$ a union of $N$ general bases, we generalize the concepts of mutual coherence and cumulative mutual coherence to the piecewise sparse case.\par
\begin{definition}
The mutual coherence of the $i$-th sub-matrix $A_i$ is
\begin{equation*}\mu^{i,i}=\max_{k\neq l}|\langle \mathbf{a}^i_k, \mathbf{a}^i_l\rangle|,~ \mathbf{a}^i_k,\mathbf{a}^i_l\in A_i,~ i=1,\ldots,N.
\end{equation*}
\end{definition}
 It is clear when $A$ is a union of orthogonal bases, $\mu_{i,i}=0$.
\begin{proposition}
The $i$-th sub-matrix coherence $\mu^{i,i}$ satisfies $0\leq \mu^{i,i}=\alpha_i \mu\leq \mu$ with $\alpha_i \in [0,1]$.
\end{proposition}
The parameter $\alpha_i$ for $i$-th block measures the ratio of coherence within $i$-th block compared with the coherence of the whole matrix.\par

\begin{definition}
The cumulative coherence between two blocks $A_i$ and $A_j$ is
\begin{equation*}
  \mu_1^{i,j}(m)=\max\limits_{|S_i|=m}~\max\limits_{l\in \{1,\ldots,n_j\}}\sum\limits_{k\in S_i}|\langle \mathbf{a}^i_k,\mathbf{a}^j_l\rangle|,
\end{equation*}
where $S_i$ is the index set of $m$ columns in sub-matrix $A_i$ and $n_j$ is the number of columns in $A_j$.
\end{definition}
\begin{remark}
Notice that the cumulative coherence between two blocks $A_i$ and $A_j$ is different from the definition of cumulative block coherence in \cite{Peotta2007}. The cumulative block coherence $\mu_{1B}(m)$ measures coherence between $m$ blocks, the $m$ represents the number of blocks.
\end{remark}
\begin{corollary}
\label{ourcoro1}
The cumulative coherence between two blocks $A_i$ and $A_j$ is bounded by
\begin{equation*}
  \mu_1^{i,j}(m)\leq m \mu.
\end{equation*}
\end{corollary}
\begin{definition}
The cumulative coherence within $A_i$ is
\begin{equation*}
  \mu_1^{i,i}(m)=\max\limits_{|S_i|=m}\max\limits_{k\notin S_i}\sum\limits_{l\in S_i}|\langle \mathbf{a}^i_k,\mathbf{a}^i_l\rangle|.
\end{equation*}
\end{definition}
\begin{remark}
Notice that $\mu_1^{i,i}(m)$ only measures the cumulative coherence within the submatrix $A_i$, ~i.e,~ how much the atoms in the same block $A_i$ are ``speaking the same language".
\end{remark}
\begin{corollary}
\label{ourcoro2}
\begin{equation*}
  \mu_1^{i,i}(m)\leq m \alpha_i \mu.
\end{equation*}
\end{corollary}

\section{Piecewise sparse recovery in union of general bases}
In the piecewise sparse setting, the P$_0$ problem (\ref{l0model}) is equivalent to the following problem:
 \begin{equation}
\label{piecewiseL0}
\begin{aligned}
\min\limits_{\mathbf{x}}&\|\mathbf{x}_1\|_0+\cdots+\|\mathbf{x}_N\|_0\\
s.t.~~&\mathbf{b}=A_1\mathbf{x}_1+\ldots+A_N\mathbf{x}_N,
\end{aligned}
\end{equation}
where $\mathbf{x}=(\mathbf{x}^T_1,\ldots,\mathbf{x}^T_N)^T$,~denote (\ref{piecewiseL0})~as piecewise P$_0$ problem.

\subsection{Uniqueness of piecewise sparse recovery via piecewise P$_0$ problem}
\begin{theorem}
\label{ourtheo1}
Suppose the measurement matrix $A=[A_1,\ldots,A_N]$ is a union of $N$ bases (or frames) with an overall coherence $\mu$ and sub--block coherence parameters $\alpha_i$ for $i=1,\ldots,N$,
if $\mathbf{x}$ is a solution of piecewise P$_0$ problem (\ref{piecewiseL0}) and
\begin{equation}\label{ourresult1}
  \|\mathbf{x}\|_0<\frac{N(1+\alpha_{\max}\mu)}{2(N-1+\alpha_{\max})\mu},
\end{equation}
then $\mathbf{x}$ is the unique solution of problem (\ref{piecewiseL0}).\par
\begin{proof}
By the sufficient and necessary condition Eq.~(\ref{sparkcondition}) for the P$_0$ problem, we need find the lower bound of the $spark(A)$ in the piecewise case.

Let $\mathbf{x}\in {\rm Ker(A)},~\mathbf{x}\neq 0$, denote the support of $\mathbf{x}$ by $S=S_1\bigcup S_2\bigcup\cdots\bigcup S_N$ and $s_i=|S_i|$, corresponding to the blocks of $A=[A_1,\ldots,A_N]$. Then
\begin{equation*}
  {\rm spark}(A)=\min\limits_{\mathbf{x}\in {\rm Ker(A)},~\mathbf{x}\neq 0}\|\mathbf{x}\|_0=\min_{\mathbf{x}\in {\rm Ker(A)},~\mathbf{x}\neq 0}(s_1+\cdots+s_N).
\end{equation*}
\par
 \textbf{Step 1.} We start similarly to the proof of Lemma 3 in \cite{gribonval2003}. Let $r_i=rank(A_i), i=1,\ldots,N.$
 Since $A_{\mathbf{S}}=[A_{S_1},\ldots,A_{S_N}]$, and
\begin{equation*}
  \mathbf{x}_{\mathbf{S}}=\left[
  \begin{array}{c}
               \mathbf{x}_{S_1} \\
               \vdots \\
               \mathbf{x}_{S_N}
             \end{array}\right]\in {\rm Ker}(A).
\end{equation*}
In order to find the minimum of $s_1+\cdots+s_N$, we can suppose that $s_i\leq r_i$ and $A_{S_i}$ is full column rank for $i=1,\ldots,N$. Because $\sum_{i=1}^N A_{S_i}\mathbf{x}_{S_i}=0$, for every $i$ we have $A_{S_i}\mathbf{x}_{S_i}=-\sum\limits_{j\neq i}A_{S_j}\mathbf{x}_{S_j}$, hence $\mathbf{x}_{S_i}=-\sum\limits_{j\neq i}(A_{S_i}^TA_{S_i})^{-1}(A_{S_i}^TA_{S_j})\mathbf{x}_{S_j}$.
Then we can deduce that
\begin{equation*}
\begin{aligned}
\|\mathbf{x}_{S_i}\|_1&\leq \frac{1}{1-\mu_1^{i,i}(s_i-1)}\sum\limits_{j\neq i}\|A_{S_i}^TA_{S_j}\|_1\|\mathbf{x}_{S_j}\|_1\\
&\leq \sum\limits_{j\neq i}\frac{\mu_1^{i,j}(s_i)}{1-\mu_1^{i,i}(s_i-1)}\|\mathbf{x}_{S_j}\|_1.
\end{aligned}
\end{equation*}
Since $\|\mathbf{x}_{\mathbf{S}}\|_1=\|\mathbf{x}_{S_1}\|_1+\cdots+\|\mathbf{x}_{S_N}\|_1$,
then
\begin{equation*}
\left(1+\frac{\max\limits_{j\neq i}\mu_1^{i,j}(s_i)}{1-\mu_1^{i,i}(s_i-1)}\right)\|\mathbf{x}_{S_i}\|_1\leq \frac{\max\limits_{j\neq i}\mu_1^{i,j}(s_i)}{1-\mu_1^{i,i}(s_i-1)}\|\mathbf{x}_{\mathbf{S}}\|_1,
\end{equation*}
which results in
\begin{equation*}
\|\mathbf{x}_{\mathbf{S}}\|_1\leq (\sum\limits_{i=1}^N\frac{v_1^i}{v_2^i})\|\mathbf{x}_{\mathbf{S}}\|_1,
\end{equation*}
where $v_1^i=\max\limits_{j\neq i}\mu_1^{i,j}(s_i)/(1-\mu_1^{i,i}(s_i-1))$ and $v_2^i=1+v_1^i$. Thus
\begin{equation}
  \sum\limits_{i=1}^N\frac{v_1^i}{v_2^i}\geq 1.
  \label{proof1}
\end{equation}
Using the inequalities:
\begin{equation*}
  \mu_1^{i,i}(s_i-1)\leq (s_i-1)\alpha_i\mu
\end{equation*}
and
\begin{equation*}
  \mu_1^{i,j}(s_i)\leq s_i\mu,
\end{equation*}
the inequality (\ref{proof1}) becomes
\begin{equation*}
  \sum\limits_{i=1}^N\frac{s_i\mu}{1-(s_i-1)\alpha_{\max}\mu+s_i\mu}\geq \sum\limits_{i=1}^N\frac{s_i\mu}{1-(s_i-1)\alpha_{i}\mu+s_i\mu}\geq 1,
\end{equation*}
where $\alpha_{\max}=\max\limits_{i=1,\ldots,N}\{\alpha_i\}$.\par
\textbf{Step 2.} In the following we evaluate the ${\rm spark}(A)$, i.e, when $s=s_1+\ldots+s_N$ reaches the minimum.\par
Denote $g_i=\frac{s_i\mu}{1-(s_i-1)\alpha_{\max}\mu+s_i\mu}$, we consider the following minimization problem
\begin{equation*}
  \min (s_1+\cdots+s_N),~~~~s.t.\sum\limits_{i=1}^Ng_i-1\geq 0.
\end{equation*}
Using the Lagrange function and KKT conditions we obtain that $s=\sum\limits_{i=1}^N s_i$ reaches minimum when $s_1=\cdots=s_N$. Then
\begin{equation*}
  \sum\limits_{i=1}^Ng_i=\frac{Ns\mu}{N(1+\alpha_{\max}\mu)+(1-\alpha_{\max})s\mu}\geq 1,
\end{equation*}
which results in
\begin{equation*}
  s\geq \frac{N(1+\alpha_{\max}\mu)}{(N-1+\alpha_{\max})\mu}.
\end{equation*}
By the definition of spark, we have $spark(A)\geq \frac{N(1+\alpha_{\max}\mu)}{(N-1+\alpha_{\max})\mu}$. Thus by Eq.~(\ref{sparkcondition}), if
\begin{equation*}
  \|\mathbf{x}\|_0<\frac{N(1+\alpha_{\max}\mu)}{2(N-1+\alpha_{\max})\mu}
\end{equation*}
then $\mathbf{x}$ is the unique solution of the piecewise P$_0$ problem (\ref{piecewiseL0}).
\end{proof}
\end{theorem}
\begin{remark}
\begin{description}
  \item[(1)] In particular, if $\alpha_{\max}=0$, i.e, $A$ is a union of $N$ orthogonal bases. The result in Theorem \ref{ourtheo1} becomes $\|\mathbf{x}\|_0<\frac{N}{2(N-1)\mu}$ which corresponds to the upper bound of Eq.~(\ref{orthouni}).\par
  \item[(2)] When~$\alpha_1=\cdots=\alpha_N=1$, i.e, $A_i$ has the same coherence with $A$. The result in Theorem \ref{ourtheo1} becomes $\|\mathbf{x}\|_0<\frac{1+\mu}{2\mu}$ which corresponds to the upper bound of Eq.~(\ref{strictbound}).\par
\end{description}
\end{remark}
\begin{example}
\label{piecewisel0example}
Consider the case when $N=2$, i.e., $A=[A_1,A_2]$. In this example we set $\mu=0.1$ and $\alpha_{\max}=0.5$, then the sufficient conditions which ensure the uniqueness for P$_0$ problem and piecewise P$_0$ problem are listed as follows:\\
\begin{enumerate}
  \item (condition 1)~$\|\mathbf{x}\|_0<\frac{1+\mu}{2\mu}$~(general condition).\\
  \item (condition 2)~$\|\mathbf{x}\|_0<\frac{1}{\mu}$~($A$ is union of orthogonal bases).\\
  \item (condition 5)~$\|\mathbf{x}\|_0<\frac{1+\alpha_{\max}\mu}{(1+\alpha_{\max})\mu}$~($A$ is union of general bases).\\
\end{enumerate}
\end{example}
\begin{figure}
  \centering
  \includegraphics[width=7.5cm]{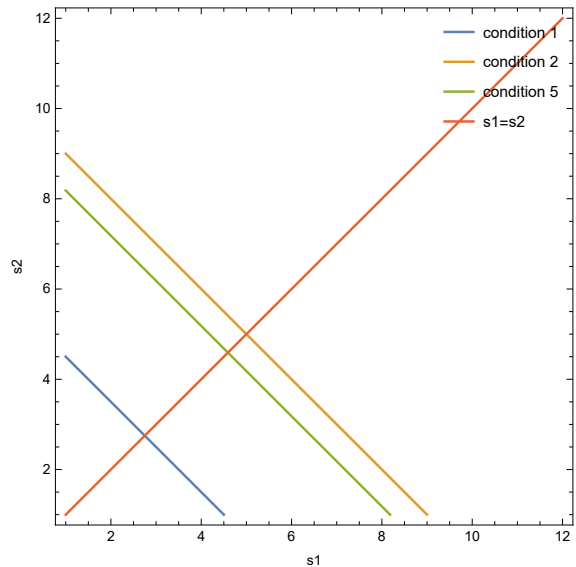}\\
  \caption{Comparison of the upper bounds for uniqueness in Example \ref{piecewisel0example}. }\label{noisefree1}
\end{figure}

From the observation of Fig.~\ref{noisefree1}, in the general case (condition 1) one can only ensure to recover $4$-sparse vector. When it comes to the piecewise sparse recovery, one can recover at least $(5,2)$-piecewise sparse vector with global $7$-sparsity by condition 5. It means that the upper bound in Theorem \ref{ourtheo1} (condition 5) is more relaxed than the upper bound in Eq.~(\ref{strictbound}) (condition 1). The condition 2 for the union of orthogonal bases is the best case. The improved condition 5
also makes a relation between general case and the union of orthogonal matrices.
Thus the results in Theorem \ref{ourtheo1} enlarge the scope the theoretical guarantees for sparse recovery by considering piecewise sparsity.\par

\subsection{Feasible conditions of algorithms for piecewise sparse recovery}
\begin{theorem}
\label{ourtheo3}
Suppose the measurement matrix $A=[A_1,\ldots,A_N]$ is a union of $N$ bases (or frames) with
an overall coherence $\mu$ and sub--block coherence parameters $\alpha_i$ for $i=1,\ldots,N$.
If
\begin{equation}\label{piecewiseERC}
2\sum\limits_{i=1}^N\frac{\mu s_i}{1+\alpha_i\mu+(1-\alpha_i)\mu s_i}< \frac{1+\alpha_Z\mu+2(1-\alpha_Z)\mu s_Z}{1+\alpha_Z\mu+(1-\alpha_Z)\mu s_Z},
\end{equation}
where $Z=\{Z: \frac{1+\alpha_Z\mu}{(1-\alpha_Z)s_Z}=\max\limits_{i=1,\ldots,N}\frac{1+\alpha_i\mu}{(1-\alpha_i)s_i}\}$, the exact recovery condition (ERC) holds.
In which case both Orthogonal Matching Pursuit and Basis Pursuit recover the sparest representation.\par
\begin{proof}
Follow the proof of Theorem 3.7 in \cite{Tropp2004} and the notations in the proof of Theorem \ref{ourtheo1}, the Grassmannian matrix
\begin{equation*}
\begin{aligned}
\Phi_{\mathbf{S}}=A_{\mathbf{S}}^{T}A_{\mathbf{S}}&=
\left(
\begin{array}{ccc}
  A_{S_{1}}^{T}  \\
   \vdots \\
    A_{S_{N}}^{T} \\
\end{array}
\right)
\left(
\begin{array}{lll}
A_{S_{1}}& \cdots & A_{S_{N}}
\end{array}
\right)\\
&=\left(
\begin{array}{cccc}
A_{S_{1}}^{T}A_{S_{1}}&A_{S_1}^TA_{S_2}&\cdots&A_{S_1}^TA_{S_N}\\
A_{S_{2}}^{T}A_{S_{1}}&A_{S_2}^TA_{S_2}&\cdots&A_{S_2}^TA_{S_N}\\
\vdots&\vdots&\ddots&\vdots\\
A_{S_{N}}^{T}A_{S_{1}}&A_{S_N}^TA_{S_2}&\cdots&A_{S_N}^TA_{S_N}\\
\end{array}
\right)=I_s-G,
\end{aligned}
\end{equation*}
where
\begin{equation*}
G=\left(
\begin{array}{cccc}
I_{s_1}-A_{S_{1}}^{T}A_{S_{1}}&-A_{S_1}^TA_{S_2}&\cdots&-A_{S_1}^TA_{S_N}\\
-A_{S_{2}}^{T}A_{S_{1}}&I_{s_2}-A_{S_2}^TA_{S_2}&\cdots&-A_{S_2}^TA_{S_N}\\
\vdots&\vdots&\ddots&\vdots\\
-A_{S_{N}}^{T}A_{S_{1}}&-A_{S_N}^TA_{S_2}&\cdots&I_{s_N}-A_{S_N}^TA_{S_N}\\
\end{array}
\right),
\end{equation*}
with the diag-block matrix $I_{s_i}-A_{S_i}^{T}A_{S_i}$ of the form $\left(
\begin{array}{cccc}
0&-A_{i_1}^TA_{i_2}&\cdots&-A_{i_1}^TA_{i_{s_i}}\\
-A_{i_2}^TA_{i_1}&0&\cdots&-A_{i_2}^TA_{i_{s_i}}\\
\vdots&\vdots&\ddots&\vdots\\
-A_{i_{s_i}}^TA_{i_1}&-A_{i_{s_i}}^TA_{i_2}&\cdots&0\\
\end{array}
\right)$.
Denote $|G|$ by the entrywise absolute value of the matrix $G$. Since all the entries in the off-diag blocks of $|G|$ can be bounded by $\mu$, and the diag-block matrix
$$|I_{s_i}-A_{S_i}^{T}A_{S_i}|\leq\left(
\begin{array}{cccc}
0&\mu^{i,i}&\cdots&\mu^{i,i}\\
\mu^{i,i}&0&\cdots&\mu^{i,i}\\
\vdots&\vdots&\ddots&\vdots\\
\mu^{i,i}&\mu^{i,i}&\cdots&0\\
\end{array}
\right),~ i=1,\cdots,N,$$
we have $|G|\leq\mu1_s-\mu B$, where $1_s$ is the $s\times s$ matrix with unit entries, $B$ is the block matrix
\begin{equation*}
B=\left(
\begin{array}{cccc}
B_1&\mathbf{0}&\cdots&\mathbf{0}\\
\mathbf{0}&B_2&\cdots&\mathbf{0}\\
\vdots&\vdots&\ddots&\vdots\\
\mathbf{0}&\mathbf{0}&\cdots&B_N\\
\end{array}
\right),
\end{equation*}
where $B_i=\alpha_iI_{s_i}+(1-\alpha_i)1_{s_i}$ is the matrix with $1$ on the diagonal, and all the off-diag entries are $1-\frac{\mu^{i,i}}{\mu}=1-\alpha_i$, $i=1,\ldots,N$.
Hence, we have the entrywise inequality
\begin{eqnarray*}
|\Phi_{\mathbf{S}}^{-1}|=|(I_s-G)^{-1}|=|I_s+\sum\limits_{k=1}^{\infty}G^k|\leq I_s+\sum\limits_{k=1}^{\infty}|G|^k\leq I_s+\sum\limits_{k=1}^{\infty}(\mu1_s-\mu B)^k\cr
=((I_s+\mu B)-\mu1_s)^{-1}=(I_s-\mu(I_s+\mu B)^{-1}1_s)^{-1}(I_s+\mu B)^{-1}.
\end{eqnarray*}
\par
\textbf{step 1}: Compute
\begin{equation*}
(I_s+\mu B)^{-1}
= \left[\begin{array}{ccc}
\frac{1}{1+\alpha_1\mu}(I_{s_1}-\frac{(1-\alpha_1)\mu}{1+\alpha_1\mu+(1-\alpha_1)\mu s_1}1_{s_1})&\cdots&\mathbf{0}\\
\vdots&\ddots&\vdots\\
\mathbf{0}&\cdots&\frac{1}{1+\alpha_N\mu}(I_{s_N}-\frac{(1-\alpha_N)\mu}{1+\alpha_N\mu+(1-\alpha_N)\mu s_N}1_{s_N})\\
\end{array}\right]
\end{equation*}

\textbf{step 2}: Compute
\begin{equation}
\label{proo1}
 (I_s-\mu(I_s+\mu B)^{-1}1_s)^{-1}=I_s+\sum\limits_{k=1}^{\infty}(\mu(I_s+\mu B)^{-1}1_s)^k,
\end{equation}
with
\begin{equation*}
\mu(I_s+\mu B)^{-1}1_s=\left[ \begin{array}{ccc} \frac{\mu}{1+\alpha_1\mu+(1-\alpha_1)\mu s_1}\mathbf{1}_{s_1}\\\vdots\\  \frac{\mu}{1+\alpha_N\mu+(1-\alpha_N)\mu s_N}\mathbf{1}_{s_N}\end{array}\right] \left[ \begin{array}{ccc}\mathbf{1}_{s_1}^T&\cdots&\mathbf{1}_{s_N}^T\end{array}\right] \overset{def}{=}\mathbf{v}\mathbf{1}_s^T.
\end{equation*}
We use $\mathbf{1}$ indicates the column vector with unit entries. Moreover, the inner product
\begin{equation*}
\mathbf{1}_s^T\mathbf{v}=\sum\limits_{i=1}^N\frac{\mu s_i}{1+\alpha_i\mu+(1-\alpha_i)\mu s_i},
\end{equation*}
therefore, the series
\begin{equation}
\sum\limits_{k=1}^{\infty}(\mathbf{v}\mathbf{1}_s^T)^k=(\mathbf{v}\mathbf{1}_s^T)\sum\limits_{k=1}^{\infty}(\mathbf{1}_s^T\mathbf{v})^{k-1}
=(\mathbf{v}\mathbf{1}_s^T)\sum\limits_{k=0}^{\infty}(\mathbf{1}_s^T\mathbf{v})^k=\frac{1}{1-\sum\limits_{i=1}^N\frac{\mu s_i}{1+\alpha_i\mu+(1-\alpha_i)\mu s_i}}(\mathbf{v}\mathbf{1}_s^T).
\end{equation}
Combined with Eq.~(\ref{proo1}), we have
\begin{equation}
\label{proo2}
|\Phi_{\mathbf{S}}^{-1}|\leq \left(I_s+\frac{1}{1-\sum\limits_{i=1}^N\frac{\mu s_i}{1+\alpha_i\mu+(1-\alpha_i)\mu s_i}}(\mathbf{v}\mathbf{1}_s^T)\right)(I_s+\mu B)^{-1}.
\end{equation}

\textbf{step 3}: Assume vector $A_i$ is drawn from basis number $Z$, then
\begin{equation*}
|A_{\mathbf{S}}^TA_j|\leq \left[ \begin{array}{ccc}|A_{S_1}^TA_j|& \cdots & |A_{S_N}^TA_j|\end{array}\right]^T \leq
\left[ \begin{array}{cccc}\mu \mathbf{1}_{s_1}^T ~ \cdots ~ \alpha_Z\mu\mathbf{1}_{s_Z}^T ~ \cdots ~ \mu \mathbf{1}_{s_N}^T \end{array}\right]^T.
\end{equation*}

So
\begin{equation}
\label{proo3}
(I_s+\mu B)^{-1}|A_{\mathbf{S}}^TA_j| \leq \left[ \begin{array}{cccc}\frac{\mu}{1+\alpha_1\mu+(1-\alpha_1)\mu s_1} \mathbf{1}_{s_1}^T ~ \cdots ~ \frac{\alpha_Z\mu}{1+\alpha_Z\mu+(1-\alpha_Z)\mu s_Z} \mathbf{1}_{s_Z}^T ~ \cdots ~ \frac{\mu}{1+\alpha_N\mu+(1-\alpha_N)\mu s_N} \mathbf{1}_{s_N}^T   \end{array}\right]^T.
\end{equation}

\textbf{step 4}: Moreover, we calculate the inner product of ERC condition $|(A_{\mathbf{S}}^TA_{\mathbf{S}})^{-1}A_{\mathbf{S}}^TA_j|$ in combination with Eqs. (\ref{proo2}), (\ref{proo3}):
\begin{equation}\label{proo4}
\begin{aligned}
|(A_{\mathbf{S}}^TA_{\mathbf{S}})^{-1}A_{\mathbf{S}}^TA_j|\leq |\Phi_{\mathbf{S}}^{-1}||A_{\mathbf{S}}^TA_j|
&\leq \left[ \begin{array}{cccc} \frac{\mu}{1+\alpha_1\mu+(1-\alpha_1)\mu s_1} \mathbf{1}_{s_1}\\ \vdots \\ \frac{\alpha_Z\mu}{1+\alpha_Z\mu+(1-\alpha_Z)\mu s_Z} \mathbf{1}_{s_Z}\\
\vdots \\ \frac{\mu}{1+\alpha_N\mu+(1-\alpha_N)\mu s_N} \mathbf{1}_{s_N} \end{array}\right]+\\
 &\frac{\sum\limits_{i\neq Z}\frac{\mu s_i}{1+\alpha_i\mu+(1-\alpha_i)\mu s_i}+\frac{\alpha_Z\mu s_Z}{1+\alpha_Z\mu+(1-\alpha_Z)\mu s_Z}}{1-\sum\limits_{i=1}^N\frac{\mu s_i}{1+\alpha_i\mu+(1-\alpha_i)\mu s_i}}
\left[ \begin{array}{ccc} \frac{\mu}{1+\alpha_1\mu+(1-\alpha_1)\mu s_1} \mathbf{1}_{s_1}\\ \vdots \\ \frac{\mu}{1+\alpha_N\mu+(1-\alpha_N)\mu s_N} \mathbf{1}_{s_N}\end{array}\right],
\end{aligned}
\end{equation}
apply the $l^1$ norm to inequality (\ref{proo4}) to reach
\begin{equation}
\label{proo5}
\|(A_{\mathbf{S}}^TA_{\mathbf{S}})^{-1}A_{\mathbf{S}}^TA_j\|_1\leq \frac{\sum\limits_{i\neq Z}\frac{\mu s_i}{1+\alpha_i\mu+(1-\alpha_i)\mu s_i}+\frac{\alpha_Z\mu s_Z}{1+\alpha_Z\mu+(1-\alpha_Z)\mu s_Z}}{1-\sum\limits_{i=1}^N\frac{\mu s_i}{1+\alpha_i\mu+(1-\alpha_i)\mu s_i}}.
\end{equation}

\textbf{step 5}: Since
\begin{equation*}
\|A_{T}^{T}A_{\mathbf{S}}^{\dagger}\|_{\infty}=\max_{j\in T}\|(A_{\mathbf{S}}^TA_{\mathbf{S}})^{-1}A_{\mathbf{S}}^TA_j\|_1,
\end{equation*}
we consider the maximum of the right side of Eq.~(\ref{proo5}), rewrite it as
\begin{equation}
\label{proo6}
\|(A_{\mathbf{S}}^TA_{\mathbf{S}})^{-1}A_{\mathbf{S}}^TA_j\|_1\leq \frac{\sum\limits_{i=1}^N\frac{\mu s_i}{1+\alpha_i\mu+(1-\alpha_i)\mu s_i}-\frac{(1-\alpha_Z)\mu s_Z}{1+\alpha_Z\mu+(1-\alpha_Z)\mu s_Z}}{1-\sum\limits_{i=1}^N\frac{\mu s_i}{1+\alpha_i\mu+(1-\alpha_i)\mu s_i}}.
\end{equation}
The right side of (\ref{proo6}) reach the maximum when $f_Z\overset{def}{=}\frac{(1-\alpha_Z)\mu s_Z}{1+\alpha_Z\mu+(1-\alpha_Z)\mu s_Z}$ reach the minimum,
\begin{eqnarray*}
\begin{aligned}
f_Z&=\frac{(1-\alpha_Z)s_Z\mu}{(1-\alpha_Z)s_Z\mu +1+\alpha_Z\mu}=\frac{\mu}{\mu+\frac{1+\alpha_Z\mu}{(1-\alpha_Z)s_Z}}.
\end{aligned}
\end{eqnarray*}
Let $Z=\{Z: \frac{1+\alpha_Z\mu}{(1-\alpha_Z)s_Z}=\max\limits_{i=1,\ldots,N}\frac{1+\alpha_i\mu}{(1-\alpha_i)s_i}\}$, and
$$2\sum\limits_{i=1}^N\frac{\mu s_i}{1+\alpha_i\mu+(1-\alpha_i)\mu s_i}<1+\frac{(1-\alpha_Z)\mu s_Z}{1+\alpha_Z\mu+(1-\alpha_Z)\mu s_Z}=\frac{1+\alpha_Z\mu+2(1-\alpha_Z)\mu s_Z}{1+\alpha_Z\mu+(1-\alpha_Z)\mu s_Z},$$
then the Exact Recovery Condition holds as $\|A_{T}^{T}A_{\mathbf{S}}^{\dagger}\|_{\infty}<1$, thus we complete the proof.
\end{proof}
\end{theorem}
In particular, when $A$ is a union of orthogonal bases, $\alpha_Z=0$, $Z$ is chosen from the minimum $s_i,~i=1,\ldots,N$, then the condition Eq.~(\ref{piecewiseERC}) corresponds to the condition Eq.~(\ref{Troppresult}) in \cite{Tropp2004}.\par
\begin{figure}
  \centering
  \includegraphics[width=7.5cm]{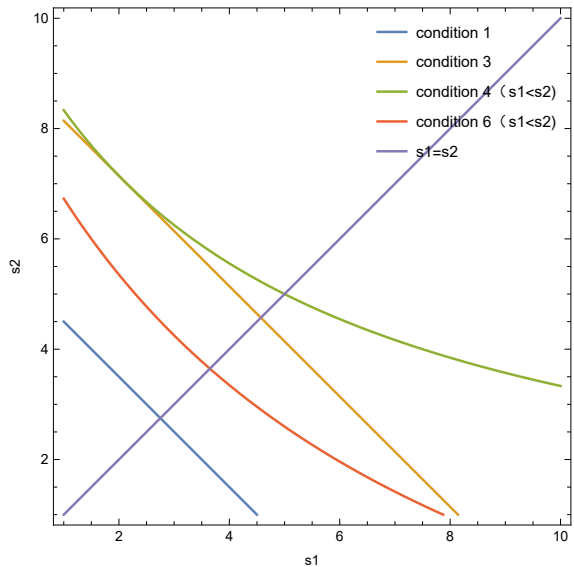}\\
  \caption{Comparison of upper bounds for feasible conditions of sparse recovery for Example \ref{ex3}.}\label{noisefree2}
\end{figure}

\begin{example}
Consider the case where $N=2$, i.e. $A=[A_1,A_2]$ and $\mathbf{x}$ is $(s_1,s_2)$-piecewise sparse vector. In this example we set overall coherence $\mu=0.1$, $\alpha_1=0.2$,~$\alpha_2=0.5$. The following sufficient conditions which ensure the feasibility of OMP and BP algorithms are listed (as shown in Fig.~\ref{noisefree2}):\\
\begin{enumerate}
  \item (condition 1)~$\|\mathbf{x}\|_0<\frac{1}{2}(1+\frac{1}{\mu})$,~(general condition)\\
  \item (condition 3)~$\|\mathbf{x}\|_0<\frac{\sqrt{2}-0.5}{\mu}$,~(equivalence condition when $A$ in pairs of orthogonal bases)\\
  \item (condition 4)
  \begin{equation*}
\begin{aligned}
&\frac{\mu s_2}{1+\mu s_2}<\frac{1}{2(1+\mu s_1)}\\
&\Leftrightarrow 2\mu^2s_1s_2+\mu s_2-1<0.
\end{aligned}
\end{equation*}
(ERC condition when $A$ in pairs of orthogonal bases)\\
  \item (condition 6)~(Assume $\alpha_1\geq \alpha_2$ and $s_1\leq s_2$, then $Z=1$)
  \begin{equation}\label{Cond6}
\begin{aligned}
&2\sum\limits_{i=1}^2\frac{\mu s_i}{1+\alpha_i\mu+(1-\alpha_i)\mu s_i}< \frac{1+\alpha_1\mu+2(1-\alpha_1)\mu s_1}{1+\alpha_1\mu+(1-\alpha_1)\mu s_1},\\
&\Leftrightarrow(1+\alpha_2\mu+(1-\alpha_2)\mu s_2)(2\alpha_1\mu s_1-\alpha_1\mu-1)+2\mu s_2(1+\alpha_1\mu+(1-\alpha_1)\mu s_1)<0.
\end{aligned}
\end{equation}
(ERC condition when $A$ in pairs of general bases).
\end{enumerate}
\label{ex3}
\end{example}

\begin{figure}
  \centering
  \includegraphics[width=8cm]{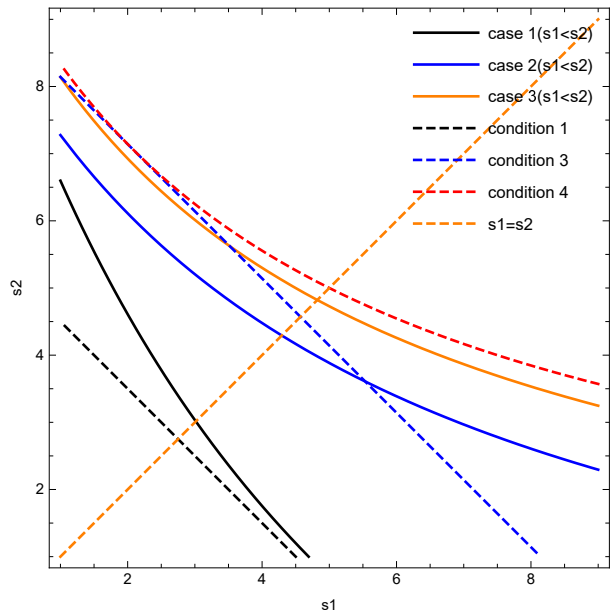}
  \caption{Comparison of upper bounds for feasible conditions of sparse recovery for Example \ref{ex4}.}\label{noisefree3}
\end{figure}

\begin{example}
In this example we show the upper bounds when $\mathbf{x}$ is $(s_1,s_2)$-piecewise sparse vector with different piecewise sparsity. We change the piecewise sparsity by vary the parameter pair $(\alpha_1,\alpha_2)$. The upper bounds of condition 6 Eq.~(\ref{Cond6}) in Example \ref{ex3} are plotted in Fig.\ref{noisefree3} for the following three cases:
\begin{enumerate}

  \item (case 1)~$\mu=0.1$,~$(\alpha_1,\alpha_2)=(0.95,0.1)$,~(the sub-matrix coherence of $A_1$ is differs greatly to that of $A_2$ and $\alpha_1$ is quite close to $1$).\\
  \item (case 2)~$\mu=0.1$,~$(\alpha_1,\alpha_2)=(0.2,0.15)$,~(the sub-matrix coherence of $A_1$ is differs slightly to that of $A_2$ and both $\alpha_1$ and $\alpha_2$ are small).\\

      \item (case 3)~$\mu=0.1$,~$(\alpha_1,\alpha_2)=(0.05,0.02)$,~(the sub-matrix coherence of $A_1$ is differs slightly to that of $A_2$ and both $\alpha_1$ and $\alpha_2$ are very small ).\\

\end{enumerate}
\label{ex4}
\end{example}

\begin{remark}
It is observed from Fig.~\ref{noisefree3} that different setting of $(\alpha_1,\alpha_2)$, i.e., different piecewise sparsity may result in different global sparsity conditions. Especially, when the sub-matrix coherences of $A_1$ and $A_2$ are small, some relaxed conditions of the sparsity can be obtained. This phenomenon provides us a guidance on the setting of piecewise structure for a given matrix in order to obtain an optimal piecewise sparsity condition, which is another interesting problem in our future work.
\end{remark}

\section{Conclusion}
In this paper, we introduce the piecewise sparsity of signals and use the mutual coherence for matrix in unions of general bases (or frames) to study the conditions for piecewise sparse recovery.
We generalize the results in orthogonal cases to the cases of general bases.
We provide the new upper bounds of global sparsity and piecewise sparsity of the signal recovered by both $l^0$ and $l^1$ optimizations when the measurement matrix $A$ is a union of general bases.
The structured information of the matrix $A$ is used to improve the sufficient conditions for successful piecewise sparse recovery and the reliability of the greedy algorithms and the BP algorithms.


%

%



\ifCLASSOPTIONcaptionsoff
  \newpage
\fi




\begin{thebibliography}{1}

\bibitem{Elad2010Book}M. Elad, \emph{Sparse and Redundant Representations: From Theory to Applications in Signal and Image Processing}, Springer, 2010.

\bibitem{mallat1} S.G.~Mallat and Z.~Zhang,``Matching pursuits with time-frequency dictionarie", \emph{IEEE Trans.
on Sig. Proces}. vol. 12, pp.3397--3415, 1993.

\bibitem{pati}Y.C. Pati, R.~Rezaiifar, and P.S.~Krishnaprasad, ``Orthogonal matching pursuit: recursive
function approximation with applications to wavelet decomposition", \emph{Proceedings of the 27th
Annual Asilomar Conference on Signals, Systems, and Computers}, 1993, pp. 40--44.

\bibitem{tropp2}J.A.Tropp and A.C. Gilbert, ``Signal recovery from random measurements via orthogonal
matching pursuit", \emph{IEEE Trans. Inform. Theory}, vol. 53, no. 12, pp. 4655--4666, 2007.

\bibitem{blumensath} T. Blumensath and M. E. Davies,``Iterative hard thresholding for compressed sensing",
\emph{Appl. Comput. Harmon. Anal.}, vol. 27, no.3, pp.265--274, 2009.

\bibitem{Donoho2012}D.L.~Donoho, Y.~Tsaig, I.~Drori, J.L.~Starck, ``Sparse solution of underdetermined linear equations by stagewise orthogonal matching pursuit (StOMP)",
\emph{IEEE Trans. on Infor. Theo.}, vol. 58, no. 2, pp.1094--1121, 2012.

\bibitem{Needell2008-1} D.Needell and R.~Vershynin, ``Uniform uncertainty principle and signal recovery via regularized orthogonal matching pursuit", \emph{Found. Comput. Math.} vol. 9, no.3, pp.317--334, 2009.

\bibitem{Needell2010} D.~Needell and R.~Vershynin, ``Signal recovery from incomplete and inaccurate measurements via regularized orthogonal matching pursuit",  \emph{IEEE J-STSP}, vol. 4, no. 2, pp.310--316, 2010.

\bibitem{Needell2008}D.~Needell and J.A~Tropp, ``CoSaMP: iterative signal recovery from incomplete and inaccurate samples",\emph{Appl. Comput. Harmon. A.}, vol. 26, no.3, pp.301--321, 2009.

\bibitem{dai}W. Dai and O. Milenkovic, ``Subspace pursuit for compressive sensing signal reconstruction",
\emph{IEEE Trans. Inform. Theory}, vol. 55, no. 5, pp.2230--2249, 2009.

\bibitem{maleki}A. Maleki, ``Coherence analysis of iterative thresholding algorithms", \emph{Proceedings of the 47th
Annual Allerton Conference on Communication, Control, and Computing}, IEEE Press, 2009,
pp. 236--243.

\bibitem{foucart}S. Foucart, ``Hard thresholding pursuit: an algorithm for compressive sensing", \emph{SIAMJ. Numer.
Anal.} vol. 49, no.6, pp. 2543--2563, 2011.

\bibitem{Kim2007}S.J. Kim, K. Koh, M. Lustig, S. Boyd, D. Gorinevsky, ``A interior--point method for large--scale $l^1$-regularized least squares", \emph{IEEE J. Select. Top. Signal Process.}, vol.1, no.4, pp. 606--617, 2007.

\bibitem{Figueiredo2007}M.A.T. Figueiredo, R.D. Nowak, S.J. Wright, ``Gradient projection for sparse reconstruction: Application to compressed sensing and other inverse problems", \emph{IEEE J. Select. Top. Signal Process.: Special Issue on Convex Optimization Methods for Signal Processing}, vol. 1, no. 4, pp. 586--598, 2007.

\bibitem{Daubechies2004}I. Daubechies, M. Defrise, C.D. Mol, ``An iterative thresholding algorithm for linear inverse problems with a sparsity constraint", \emph{Comm. Pure Appl. Math.}, vol. 57, pp. 1413--1457, 2004.

\bibitem{Donoho2001}D.L. Donoho and X. Huo, ``Uncertainty principles and ideal atomic decomposition",
\emph{ IEEE Trans. It.}, vol. 47, no.7, pp. 2845-2862, 2001.

\bibitem{Donoho2003}D.L. Donoho and M. Elad, ``Optimally sparse representation in general (nonorthogonal) dictionaries via $l^1$ minimization", \emph{Proceedings of the National Academy of Sciences of the United States of America}, vol. 100, no.5, 2003, pp. 2197-2202.

\bibitem{Tropp2004}J.A Tropp, ``Greed is good: algorithmic results for sparse approximation", \emph{IEEE T. Inform. Theory}, vol. 50, no. 10, pp. 2231--2242, 2004.

\bibitem{Candes2005}E.J. Canddes and T. Tao, ``Decoding by linear programming". \emph{IEEE T. Inform. Theory},  vol. 51, pp. 4203-4215, 2005.

\bibitem{Candes2006-2}E.J. Canddes, J.K. Romberg, T. Tao, ``Stable signal recovery from incomplete and inaccurate measurements", \emph{Communications on Pure \& Applied Mathematics}, vol. 59, no.8, pp.1207--1223, 2006.

\bibitem{Candes2008-1}E.J. Canddes, ``The restricted isometry property and its implications for compressed sensing", \emph{Comptes rendus-Mathsmatique}, vol. 346, no. 9, pp. 589--592, 2008.

\bibitem{elad2002}M. Elad and A. M. Bruckstein,``A generalized uncertainty principle and sparse
representation in pairs of bases", \emph{IEEE T. Inform. Theory},  vol. 48, no. 9, pp. 2558--2567, 2002.

\bibitem{gribonval2003} R. Gribonval and M. Nielsen, ``Sparse representations in unions of bases", \emph{IEEE T. Inform. Theory}, vol. 49, no.12, pp. 3320--3325, 2003.

\bibitem{Peotta2007}L.~Peotta and P.~Vandergheynst,``Matching pursuit with block incoherent dictionaries", \emph{IEEE Trans. Signal Proces.}, vol. 55, no. 9, pp. 4549--4557, 2007.

\bibitem{Eldar2009}Y.C.~Eldar and M,~ Mishali, ``Block sparsity and sampling over a union of subspaces",in \emph{International Conference on Digital Signal Processing}. IEEE, 2009, pp. 1--8.

\bibitem{Eldar2010} Y.C.~Eldar, P.~Kuppinger, H.~Bolcskei, `` Block-sparse signals: uncertainty relations and efficient recovery", \emph{IEEE Trans. Signal Proces.},vol. 58, no.6, pp.3042--3054, 2010.

\bibitem{Elhamifar2011} E.~Elhamifar and R.~Vidal, ``Block-sparse recovery via convex optimization", \emph{IEEE Trans.Signal Proces.},vol. 60, no. 8, pp. 4094--4107, 2011.


\bibitem{Starck2005}J.L. Starck, M. Elad and D.L. Donoho, ``Image decomposition via the combination of sparse representations and a variational approach", \emph{IEEE Trans. image proces.},vol. 14, no. 10, pp. 1570--1582, 2010.

    \bibitem{hlw}Y.X. Hao, C.J. Li and R.H. Wang, ``Sparse approximate solution of fitting surface to scattered points by MLASSO model",\emph{Sci. China Math.},  vol. 7, pp. 1--18, 2018.

        \bibitem{Zhong2018} Y.J. Zhong, C.J. Li, ``Piecewise sparse recovery via piecewise inverse
scale space algorithm with deletion rule", Journal of Computational Mathematics, to be published.





\end{thebibliography}
%

\end{document}